\documentclass[11pt]{article}  
\usepackage{amsfonts}
\usepackage{moresize}
\usepackage{hyperref}
\usepackage{mathtools}
\usepackage[top=3cm]{geometry}
\usepackage{tikz}
\begin{document}

\newcommand{\li}{\underline}
\newcommand{\PP}{{\cal P}}
\newcommand{\RR}{{\cal R}}
\newcommand{\lam}{\lambda}
\newcommand{\ei}{\epsilon_i}
\newcommand{\be}{\bar{\epsilon}}
\newcommand{\eproof}{\hfill $\Box$ \vspace{0.2cm}}
\newcommand{\qed}{\hfill$\square$}
\newcommand{\bra}[1]{\langle#1|} 
\newcommand{\ket}[1]{|#1\rangle} 
\newcommand{\braket}[2]{ \langle #1 | #2 \rangle} 

\newtheorem{theorem}{Theorem}
\newtheorem{proposition}{Proposition}
\newtheorem{rem}{Remark}
\newtheorem{corol}{Corollary}
\newtheorem{defi}{Definition}
\newtheorem{notat}{Notation}
\newtheorem{lemm}{Lemma}
\newtheorem{exam}{Example}

\vspace{-\baselineskip}
\vspace{-\baselineskip}
\vspace{-\baselineskip}
\vspace{-\baselineskip}
\vspace{-\baselineskip}
\vspace{-\baselineskip}

\title{\bf Conformal dimensionality reduction / increase\\
\Large{\textcolor{gray}{(up to an arbitrarily chosen dimension number)}}}

\author{Nicholas J. Daras\\
\small{\textit{Department of Mathematics, Hellenic Army Academy, 16673 Attiki, Greece}} \\
\small{Email : njdaras@sse.gr }}
\date{}
\maketitle
\setlength{\parindent}{1pt}
\vspace{-\baselineskip}

\begin{footnotesize}
\noindent \textbf{Abstract}. \textit{We give two low-complexity algorithms, one for dimensionality reduction and one for dimensionality increase, which are applicable to any dataset, regardless of whether the set has an intrinsic dimension or not. The corresponding methods introduce chains of compositions of conformal homeomorphisms that transform any data set $\mathbb{X}$ in a Euclidean space  $\mathbb{R}^{D+1}$ into an isopleth dataset $ \mathbb{Y}$ within a Euclidean space $\mathbb{R}^{\mathfrak{D}+1}$ of arbitrarily smaller or of arbitrarily larger dimension $\mathfrak{D}+1$ and preserve all angles, in the sense that all angles formed between points in the original dataset  $ \mathbb{X}$ are equal to the angles formed between the images of these points in the new dataset  $\mathbb{Y}$. Because they preserve angles, the two methods also preserve  shapes locally, although, in general, the overall sizes and shapes are distorted away from a center point.}  \\


\end{footnotesize}

\setlength{\parindent}{6ex} {\Large \textbf{Introduction}}\\

\setlength{\parindent}{6ex} The (non-linear) \textbf{\textit{dimensionality reduction problem}} can be defined as follows (\cite{CarreiraPerpinan}, \cite{Fodor},\cite{ VanDerMateenPostmaVanDenHerik}). Assume we have a dataset
\begin{center}
$\mathbb{X}_{D+1}=\lbrace \textbf{x}^{(\nu )}=\left( x_{1}^{(\nu )}, x_{2}^{(\nu )},...,x_{D}^{(\nu )},x_{D+1}^{(\nu )} \right) \in \mathbb{R}^{D+1}, \nu = 1, 2,..., n \rbrace \left(\subset \mathbb{R}^{D+1} \right)$ 
\end{center}
consisting of $n$ data vectors $ \textbf{x}^{(\nu )} $ ($\nu \in \lbrace 1, 2,..., n \rbrace$) all of which lie within a Euclidean space $\mathbb{R}^{D+1}$ of dimension $D+1$. Assume further that this dataset has \textit{intrinsic dimension}  \footnote {The intrinsic dimension for a dataset can be thought of as the minimal number of variables needed to represent the dataset. Similarly, in signal processing of multidimensional signals, the intrinsic dimension of the signal describes how many variables are needed to generate a good approximation of the signal. The definition of intrinsic dimension is as follows. Let $f=f(\textbf{x})$ be a $D-$variable function (:the set of variables can be represented as a $D-$dimensional vector $\textbf{\textbf{x}}$). If for some $d-$variable function $g$ and $d\times D $ matrix $M$ it is the case that
\begin{itemize}
\item $f(\textbf{x})=g(M\textbf{x})$ for all $\textbf{x}$ and 
\item $ d$ is the smallest number for which the above relation between $f$ and $g$ can be found,
\end{itemize}
then \textit{the intrinsic dimension of} $f$ is $d$. For instance, let $f\left(x_{1}, x_{2} \right) = g \left(x_{1}+ x_{2} \right)$. $f$ is still intrinsic one-dimensional, which can be seen by making a variable transformation $y_{1}=x_{1}+ x_{2} $ and $y_{2}=x_{1}- x_{2} $ which gives 
\begin{center}
$f\left( \frac{y_{1}+y_{2}}{2},\frac{y_{1}-y_{2}}{2}\right)=g\left( y_{1} \right) $.
\end{center}
Since the variation in $f$ can be described by the single variable $y_{1}$ its intrinsic dimension is one.
}
$d+1$ (with $d<D$). Here, in mathematical terms, intrinsic dimension means that the points in dataset $\mathbb{X}_{D+1}$ are lying on or near a manifold with dimension $d+1$ that is embedded in the $(D+1)-$dimensional space.\\

 Note that we make no assumptions on the structure of this manifold: the manifold may be non-Riemannian because of discontinuities (i.e., the manifold may consist of a number of disconnected submanifolds). \\
\indent \textit{Dimensionality reduction techniques transform the dataset} $\mathbb{X}_{D+1}\subset \left( \mathbb{R}^{D+1} \right) $  \textit{into a new dataset }
\begin{center}
$\mathbb{Y}_{d+1}=\lbrace \textbf{y}^{(\nu )}=\left( y_{1}^{(\nu )}, y_{2}^{(\nu )},...,y_{d}^{(\nu )},y_{d+1}^{(\nu )} \right) \in \mathbb{R}^{d+1}, \nu = 1, 2,..., n \rbrace \left(\subset \mathbb{R}^{d+1} \right)$ 
\end{center} 
\textit{with} ($1\leq$ ) $ d<D$, while \textit{retaining the geometry of the data as much as possible}. \\
\indent However, in general, \textit{neither the geometry of the data manifold, nor the intrinsic dimension $d+1$ of the dataset $\mathbb{X}_{D+1}$ are known}. Therefore, \textit{dimensionality reduction seems to be an ill-posed problem} that can only be solved by assuming certain properties of the data (such as its intrinsic dimension). In fact, and in order to be consistent with the general case and adequately prepared for all possible situations that may arise in practice, we must clarify that \textit{intrinsic dimension does not always exist}   or, even if it does exist, \textit{is not always easily estimated for every dataset or system} \footnote{Here are some indicative reasons why the intrinsic dimension may not always exist or not be easily estimated.
\begin{itemize}
\item \textit{Scale-Dependence}: The estimated intrinsic dimension can vary significantly depending on the scale at which you examine the data.
\item \textit{Noise}: Noise in the data can lead to overestimation of the intrinsic dimension, especially at short distances between points.
\item \textit{Complex Geometries}: Some geometric objects, like space-filling curves, pose formal ambiguities, making a universally consistent definition of intrinsic dimension difficult.
\item \textit{Varying Definitions}: Different definitions of intrinsic dimension exist, and the appropriate one depends on the specific problem, leading to a ``spectrum of various dimensions" rather than a single, universal concept.
\item \textit{Local vs. Global}: The intrinsic dimension might vary locally across a dataset, with some methods handling local dimensions (:local intrinsic dimension) while others aim for a global measure.
\end{itemize}
}. While it's a useful concept for understanding the complexity of high-dimensional data, especially when data lies on a lower-dimensional manifold, the ``existence" and value of intrinsic dimension can depend on factors like the presence of noise, the scale of the data, and the formal definition used, with some theoretical challenges like space-filling curves making a universal definition problematic. In Section 2, we will completely avoid the concept of intrinsic dimension and give a new \textit{low-complexity dimensionality reduction algorithm} that is applicable in every case, without being at all concerned with ensuring conditions under which intrinsic dimension exists. More specifically, the implementation of our method introduces a \textit{chain of compositions of conformal homeomorphisms} that 
\begin{itemize}
\item \textit{transforms any data set $\mathbb{X}_{D+1}$ $(\subset\mathbb{R}_{D+1})$  within a Euclidean space $\mathbb{R}_{D+1}$ of dimension $D+1$ into a new (isopleth) dataset $\mathbb{Y}_{d+1}$ $(\subset\mathbb{R}_{d+1})$ within another Euclidean space $\mathbb{R}_{d+1}$ of arbitrarily chosen smaller dimension} $d+1$ (: $D > d \geq 1$) and 
\item \textit{preserves angles}, in the sense that \textit{all angles formed between the points of the original dataset $\mathbb{X}_{D+1}$ are equal to the angles formed between the images of these points within the new dataset} $\mathbb{Y}_{d+1}$. 
\end{itemize}
Because it preserves angles,
\begin{itemize}
\item \textit{this method also locally preserves shapes}, although the overall size and  shape of the original dataset are distorted away from a center point.
\end{itemize}

\begin{rem}
\normalfont{Throughout the paper,} we do not assume the dataset $\mathbb{X}_{D+1}$ is zero-mean \normalfont{ (contrary to the usual assumption; see, for instance, \cite{ VanDerMateenPostmaVanDenHerik}}). $\square$ 
\end{rem}

Although in many machine learning methods, we usually try to reduce the dimensionality and find a latent space/manifold in which the data can be represented (like in neural networks taking images), in other methods like SVM/kernel \footnote{SVM's are kernel methods. A kernel can be in an infinite dimensional feature space. After evaluating data points in this ``higher" dimension, we can then perform linear inference on them (using Mercer’s Theorem). SVM's specifically use the higher dimensionality to aid in drawing a hyperplane which can separate the data (Kernel trick).}, we are trying to find a higher dimensional space so we can separate/classify our data.\\
\indent More generally, in data science, one can argue that increasing dimensionality has some notable benefits, such as
\begin{enumerate}
\item \textit{improved data representation} : the fact that higher dimensions allow for a more detailed and complex representation of data, capturing intricate patterns and relationships that may not be observed in lower-dimensional spaces,
\item \textit{improved model accuracy}  \footnote{up to a point} : adding relevant features can provide additional information that improves classifier performance, although this benefit diminishes and can be reversed beyond an optimal number of features \footnote{\textit{Hughes Effect}}  and
\item at its core, vector dimensionality refers to the number of components or features that make up a vector, so, as the dimensionality increases, a vector can represent more complex data structures, allowing for \textit{more nuanced similarities}.
\end{enumerate}

Particularly in machine learning problems that involve learning a ``natural state" from a finite number of data samples in a high-dimensional feature space, with each feature having a range of possible values, \textit{a huge amount of training data is usually required} to ensure that there are many samples with every combination of values.\\

Nevertheless, various phenomena that arise when analysing and organizing data in high-dimensional spaces that do not occur in low-dimensional environments and which may create computational problems. The common theme of these problems is that when the dimensionality increases, \textit{the volume of the space increases so fast that the available data become sparse}. In order to obtain a reliable result, the amount of data needed often grows exponentially with the dimensionality. Also, organizing and searching data often relies on detecting areas where objects form groups with similar properties; in high dimensional data, however, all objects appear to be sparse and dissimilar in many ways, which prevents common data organization strategies from being efficient. Briefly, one can say that the challenges of dimensionality increase (: \textit{the curse of dimensionality}; see Richard Ernest Bellman's book \cite{Bellman}) ) are:

\begin{enumerate}
\item \textit{data sparsity} : as dimensionality increases, the volume of the space expands exponentially, making data points much further apart and sparse; this means that even large datasets can become ``thinly" populated, making it difficult for algorithms to find meaningful patterns,
\item \textit{distance metrics lose meaning} : algorithms that rely on distances, such as k-nearest neighbours (k-NN), find it hard to differentiate between ``near" and ``far" points in high-dimensional spaces because the distances become more uniform,
\item \textit{model over-fitting} : with more dimensions, models require more parameters, increasing the risk of overfitting to noise in the data rather than capturing the true underlying structure, especially if the dataset is not sufficiently large and
\item  \textit{increased computational cost} : higher dimensions require more memory and processing power, leading to slower training times and higher computational requirements.
\end{enumerate}

The \textbf{\textit{dimensionality increase problem}} can be defined analogously to the dimensionality reduction problem: \textit{Transform a given dataset }
\begin{center}
$\mathbb{X}_{D+1}=\lbrace \textbf{x}^{(\nu )}=\left( x_{1}^{(\nu )}, x_{2}^{(\nu )},...,x_{D}^{(\nu )},x_{D+1}^{(\nu )} \right) \in \mathbb{R}^{D+1}, \nu = 1, 2,..., n \rbrace \left(\subset \mathbb{R}^{D+1} \right)$
\end{center}
\textit{into a new dataset} 
\begin{center}
 $\mathbb{Y}_{D^{'}+1}=\lbrace \textbf{y}^{(\nu )}=\left( y_{1}^{(\nu )}, y_{2}^{(\nu )},...,y_{D^{'}}^{(\nu )},y_{D^{'}+1}^{(\nu )} \right) \in \mathbb{R}^{D^{'}+1}, \nu = 1, 2,..., n \rbrace \left(\subset \mathbb{R}^{D^{'}+1} \right)$
\end{center} 
\textit{with}  $D^{'}>D$, while \textit{retaining the geometry of the data as much as possible}. In Section 3, we will give a low-complexity dimensionality increase algorithm that is applicable in every case. Its implementation introduces a chain of compositions of conformal homeomorphisms that 
\begin{itemize}
\item \textit{transforms any dataset $\mathbb{X}_{D+1} $ within a Euclidean space $\mathbb{R}^{D+1}$ of dimension $D+1$ into a new isopleth dataset $\mathbb{Y}_{D^{'}+1} $ within another Euclidean space $\mathbb{R}^{D^{'}+1}$ of arbitrarily chosen larger dimension} $D^{'}+1 (>D+1)$ and 
\item \textit{preserves angles, in the sense that all angles formed between the points of the original dataset $\mathbb{X}_{D+1}$ are equal to the angles formed between the images of these points within the new dataset} $\mathbb{Y}_{D^{'}+1}$. 
\end{itemize}
Because it preserves angles,
\begin{itemize}
 \item \textit{our method also locally preserves shapes}, although the overall size and  shape of the original dataset are distorted away from a center point.
\end{itemize} 

The key-concept for the formulation and construction of the two algorithms given in Sections 2 and 3 is the \textit{generalized stereographic projection} from the (North Pole-punctured) hypersphere of dimension $D>3$ to the corresponding Euclidean space of dimension $D+1$, as well as the definition of \textit{inverse generalized stereographic projection} from the Euclidean space of dimension $D+1$ ($D\geq 1$) to the (North Pole-punctured) hypersphere of dimension $D$. In Section 1 we will introduce these two mappings in a reasonable way. As we will see, the two fundamental properties of both are
\begin{itemize}
\item on the one hand, that \textit{they are homemorphisms} and,
\item on the other hand, that \textit{they preserve the angles formed between any points of the domains of definition in the angles formed between the images of these points}.
\end{itemize}
Based on the constructive definition of these two mappings and exploiting their two aforementioned fundamental properties we will arrive at the formulation of the two algorithms in Sections 2 and 3.

\section{Generalized stereographic projection}

Let 
\begin{center}
$\mathcal{S}^{D}=\lbrace (\chi_{1}, \chi_{2},..., \chi_{D}, \chi_{D+1})\in \mathbb{R}^{D+1}: \sum_{\nu=1}^{D+1}\chi_{\nu}^{2}=1 \rbrace$ 
\end{center}
be the unit hypersphere (also called $D-$sphere) of the $D+1-$dimensional Euclidean space $\mathbb{R}^{D+1}$. If $\mathsf{N}:=(0,0,...0,1)\in \mathbb{R}^{D+1}$ is its \textit{North Pole} and if $\mathbb{R}^{D}$ is taken as a hyperplane (: $D-$plane) in $\mathbb{R}^{D+1}$, then the \textit{generalized stereographic projection} of a point $\Sigma \in \mathcal{S}^{D} \setminus \lbrace \mathsf{N} \rbrace$ is the point $\Sigma^{'}$ of intersection of the line $N\Sigma$ with the hyperplane $\mathbb{R}^{D}$. In Cartesian coordinates, the \textit{generalized stereographic projection} $\mathcal{P}_{D}$ from $\mathcal{S}^{D} \setminus \lbrace \mathsf{N} \rbrace$ onto $\mathbb{R}^{D}$  is given by the mapping \\

$\mathcal{P}_{D}:\mathcal{S}^{D} \setminus \lbrace \mathsf{N} \rbrace \longrightarrow \mathbb{R}^{D}$\\

$ (\chi_{1}, \chi_{2},..., \chi_{D}, \chi_{D+1}) \mapsto \mathcal{P}_{D}(\chi_{1}, \chi_{2},..., \chi_{D}, \chi_{D+1})
:= 
\biggl(
\frac{\chi_{1}}{\underbrace{1-\chi_{D+1}}_{:=x_{1}}},\frac{\chi_{2}}{\underbrace{1-\chi_{D+1}}_{:=x_{2}} },...,
\frac{\chi_{D}}{\underbrace{1-\chi_{D+1}}_{:=x_{D}}}
\biggr).
$\\

\noindent Indeed, if $\Sigma_{0}^{'}$ is the point where the line through $\mathsf{N}$ to the point $\Sigma \in \mathcal{S}^{D} \setminus \lbrace \mathsf{N} \rbrace$ intersects the hyperplane $\mathbb{R}^{D}=\lbrace x_{D+1}= 0\rbrace$ in $\mathbb{R}^{D+1}$, then the point $\Sigma_{0}^{'}$ is characterized by the fact that the straight line segment $\overline{\mathsf{N}\Sigma_{0}^{'}}$ equals a real multiple (dilation or contraction) of $\overline{\mathsf{N}\Sigma}$:
\begin{center}
 $\overline{\mathsf{N}\Sigma_{0}^{'}} = \lambda \overline{\mathsf{N}\Sigma}$
\end{center}
for some non-zero real scalar $\lambda$. Writing $\Sigma_{0}^{'}=(x_{1},x_{2},…,x_{D},0)$ and $\Sigma= (\chi_{1}, \chi_{2},..., \chi_{D}, \chi_{D+1})\in \mathbb{R}^{D+1}$, this leads to the system of equations
\begin{equation}
\left.\begin{aligned}
x_{\nu}  &=\lambda \chi_{\nu} \hspace{0.1cm}(\nu=1,2,...,D),\\
   1     &=\lambda (1-\chi_{D+1}).
        \end{aligned}\\
\right\}
\qquad \text{}
\end{equation}
Solving the second equation for $\lambda (\neq 0)$ and plugging it into the first equation, we obtain the above formula for the stereographic projection, that is
\begin{equation}
\mathcal{P}_{D}(\chi_{1}, \chi_{2},..., \chi_{D}, \chi_{D+1})
=
\biggl(
\frac{\chi_{1}}{\underbrace{1-\chi_{D+1}}_{:=x_{1}}},\frac{\chi_{2}}{\underbrace{1-\chi_{D+1}}_{:=x_{2}} },...,
\frac{\chi_{D}}{\underbrace{1-\chi_{D+1}}_{:=x_{D}}}
\biggr).
\end{equation}

To compute the inverse $\mathcal{P}_{D}^{-1}$ of $\mathcal{P}_{D}$, we may eliminate $\lambda$ in (1) and obtain
\begin{equation}
x_{\nu}=\frac{\chi_{\nu}}{1-\chi_{D+1}} \hspace{0.1cm}(\nu=1,2,...,D);
\end{equation}
But the above relation only gives us $D$ equations and we have $(D+1)$ unknowns. To find one more equation, we can do some reverse-engineering. For $\mathcal{P}_{D}^{-1}$ to be well defined, it is necessary that the output be a point of the sphere. This gives us the extra equation we needed:
\begin{equation}
\chi_{1}^{2}+\chi_{2}^{2}+...+\chi_{D}^{2}+\chi_{D+1}^{2}=1.
\end{equation}
Plugging (1) in (4), we take $1=(1-\chi_{D+1})^{2} (x_{1}^{2}+x_{2}^{2}+...+x_{D}^{2})+\chi_{D+1}^{2}=(1-\chi_{D+1})^{2} \mid\mid \textbf{x}\mid\mid^{2}+\chi_{D+1}^{2}$, where $\mid\mid \textbf{x}\mid\mid^{2}:=\sum_{\nu=1}^{D} x_{\nu}^{2}=x_{1}^{2}+x_{2}^{2}+...+x_{D}^{2}$. This is a second order equation with $\chi_{D+1}$ being the only unknown:
\begin{center}
$(1-\mid\mid \textbf{x}\mid\mid^{2}) \chi_{D+1}^{2}-2\mid\mid \textbf{x}\mid\mid^{2}\chi_{D+1}+(\mid\mid \textbf{x}\mid\mid^{2}-1)=0$.
\end{center}
Applying the classical formula to solve this equation, you either get $\chi_{D+1}=1$ or $\chi_{D+1}=1-2/((1+\mid\mid \textbf{x}\mid\mid^{2})$. Since we know that must satisfy the initial relations, it is obvious that the first solution must be discarded, so
\begin{center}
$\chi_{D+1}=1-\frac{2}{1+\mid\mid \textbf{x}\mid\mid^{2}}$.
\end{center}
Plugging this in the equations (2), we get:
\begin{equation}
\left.\begin{aligned}
\chi_{\nu}  &=\frac{2x_{\nu}}{\mid\mid \textbf{x}\mid\mid^{2}} \hspace{0.1cm}(\nu=1,2,...,D),\\
\chi_{D+1}  &=1-\frac{2}{1+\mid\mid \textbf{x}\mid\mid^{2}}=\frac{\mid\mid \textbf{x}\mid\mid^{2}-1}{\mid\mid \textbf{x}\mid\mid^{2}+1}.
        \end{aligned}\\
\right\}
\qquad \text{}
\end{equation}
Inserting (5) back into (3) gives the formula
\begin{center}
$\mathcal{P}_{D}^{-1}(x_{1},x_{2},..,x_{D})
=
\biggl(
\frac{2 x_{1}}{\underbrace{\mid\mid \textbf{x}\mid\mid^{2}+1}_{:=\chi_{1}}},
\frac{2 x_{2}}{\underbrace{\mid\mid \textbf{x}\mid\mid^{2}+1}_{:=\chi_{2}}},...,
\frac{2 x_{D}}{\underbrace{\mid\mid \textbf{x}\mid\mid^{2}+1}_{:=\chi_{D}}},
\frac{\mid\mid \textbf{x}\mid\mid^{2}-1}{\underbrace{\mid\mid \textbf{x}\mid\mid^{2}+1}_{:=\chi_{D+1}}}
\biggr)
$
\end{center}
which maps $\mathbb{R}^{D}$ back onto $\mathcal{S}^{D}\setminus{N}$.\\
 
Since $\mathcal{P}_{D}$ and $\mathcal{P}_{D}^{-1}$ are both visibly compositions of continuous mappings, we have the

\begin{theorem}
The generalized stereographic projection of the punctured hypersphere $\mathcal{S}^{D} \setminus \lbrace \mathsf{N} \rbrace$ onto $\mathbb{R}^{D}$ is the map\\

$\mathcal{P}_{D}: (\mathbb{R}^{D+1}\supset)\hspace{0.1cm} \mathcal{S}^{D} \setminus \lbrace \mathsf{N} \rbrace \longrightarrow \mathbb{R}^{D} $\\

$ (\chi_{1}, \chi_{2},..., \chi_{D}, \chi_{D+1}) \mapsto \mathcal{P}_{D}(\chi_{1}, \chi_{2},..., \chi_{D}, \chi_{D+1})
:= 
\biggl(
\frac{\chi_{1}}{\underbrace{1-\chi_{D+1}}_{:=x_{1}}},\frac{\chi_{2}}{\underbrace{1-\chi_{D+1}}_{:=x_{2}} },...,
\frac{\chi_{D}}{\underbrace{1-\chi_{D+1}}_{:=x_{D}}}
\biggr).
$\\

\noindent It is a homeomorphism with inverse from $\mathbb{R}^{D}$ onto $\mathcal{S}^{D} \setminus \lbrace \mathsf{N} \rbrace$ given by the map:\\

$\mathcal{P}_{D}^{-1}:\mathbb{R}^{D}\longrightarrow \mathcal{S}^{D} \setminus \lbrace \mathsf{N} \rbrace \hspace{0.1cm} (\subset \mathbb{R}^{D+1})$\\

$ (x_{1},x_{2},..,x_{D}) \mapsto \mathcal{P}_{D}^{-1}(x_{1},x_{2},..,x_{D})
=
\biggl(
\frac{2 x_{1}}{\underbrace{\mid\mid \textbf{x}\mid\mid^{2}+1}_{:=\chi_{1}}},
\frac{2 x_{2}}{\underbrace{\mid\mid \textbf{x}\mid\mid^{2}+1}_{:=\chi_{2}}},...,
\frac{2 x_{D}}{\underbrace{\mid\mid \textbf{x}\mid\mid^{2}+1}_{:=\chi_{D}}},
\frac{\mid\mid \textbf{x}\mid\mid^{2}-1}{\underbrace{\mid\mid \textbf{x}\mid\mid^{2}+1}_{:=\chi_{D+1}}}
\biggr),
$
where we have set
\begin{center}
$\mid\mid \textbf{x}\mid\mid^{2}:=\sum_{\nu=1}^{D} x_{\nu}^{2}=x_{1}^{2}+x_{2}^{2}+...+x_{D}^{2}$. $\square $\\
\end{center}
\end{theorem}

\vspace{5mm}

Generalized stereographic projection's key properties are its ability to map a $D-$sphere to a $D-$plane in a \textit{conformal} (:\textit{angle-preserving})\footnote{A \textit{conformal mapping} in $D-$dimensional Euclidean space is a continuous transformation that preserves the angles between any two intersecting curves at their point of intersection, both in magnitude and orientation. This angle-preserving property also implies a constant local scaling factor, meaning that the ratio of the distances between a mapped point and a reference point to the original distance tends to a fixed limit.} and \textit{bijective} (:\textit{one-to-one and onto}) manner, mapping circles to circles or lines, and preserving angles locally while distorting area and distance globally. It is a fundamental tool in complex analysis for creating the Riemann sphere and in differential geometry for parametrizing the hypersphere. \\

\newpage

\begin{theorem}
\normalfont{\textbf{Generalized stereographic projection's key properties}:}
\begin{enumerate}
\item \normalfont{\textbf{Conformal Mapping}\footnote{The three characteristic properties of a conformal mapping are the following.
\begin{itemize}
\item \textbf{Angle Preservation}: The most fundamental aspect is the conservation of angles between curves. If two curves meet at a certain angle, their images under the conformal map will also meet at the same angle at the corresponding point.
\item \textbf{Dilation/Scaling}: Conformal maps also exhibit a constant local scaling factor. The ratio of the distance between mapped points to the original distance approaches a constant value, known as the coefficient of dilation.
\item \textbf{Local Transformation}: A map is considered conformal if this angle-preserving property holds at every point within a given domain.
\end{itemize}
Conformal mappings are widely used in fields like complex analysis, fluid dynamics, and cartography (e.g., the \textit{Mercator projection}) to simplify complex geometries and solve problems in two-dimensional spaces. Notice that in the context of complex numbers, \textit{a function $f(z)$ is conformal if it is analytic and its derivative $f'(z)$ is non-zero for any $z$ in the domain}. This mathematical property ensures the \textit{preservation of angles} and is a cornerstone of the theory of complex functions.}:} \textit{The most significant property of generalized stereographic projection is that it is a conformal mapping} (see also Theorem 3 below). 
\item \textbf{Preserves Angles} \footnote{In Euclidean space, an \textit{angle} is fundamentally the geometric figure formed by two rays (or lines) sharing a common endpoint called the vertex. The measure of this ``opening" is defined by the dot product of the two vectors representing these rays, using the formula $\theta = Arccos((u \cdot v) / (\parallel u \parallel \parallel v \parallel))$, where $u$ and $v$ are the vectors, and $\theta$ is the angle. This relationship allows angles to be quantified using units like degrees or radians, which represent a fraction of a full circle or rotation.}:  \textit{It preserves the angles at which curves intersect. If two curves on the hypersphere intersect at an angle $\theta$, their images on the hyperplane will also intersect at the same angle} (see also Theorem 3 below). 
\item \textbf{Preserves Shapes} (Locally): \textit{Because it preserves angles, it locally preserves shapes, though the overall size and shape are distorted away from the center point}.
\item \textbf{Bijective Mapping}: \textit{The projection creates a one-to-one correspondence between the hypersphere and the hyperplane}. 
\item \textbf{Maps Sphere to Plane}: \textit{It maps the entire hypersphere, except for the point of projection (the ``North Pole"), onto the entire hyperplane}. 
\item \textbf{The Point at Infinity}: \textit{The point of projection on the hypersphere corresponds to the point at infinity on the hyperplane}. 
\item \textbf{Circles and Lines}: \textit{A specific geometric property is how it handles circles}. 
\item \textbf{Circles to Circles}: \textit{Circles on the hypersphere that do not pass through the center of projection map to circles on the hyperplane}.
\item \textbf{Circles Through Pole to Lines}: \textit{Circles on the hypersphere that pass through the center of projection map to straight lines on the hyperplane}. $\square $
\end{enumerate}
\end{theorem}

\begin{rem}
The generalized stereographic projection is neither isometric nor equiareal. \normalfont{Indeed, generalized stereographic projection does not preserve distances (: distances between points on the sphere are distorted on the plane).  Similarly, it is not an equiareal projection, meaning it does not preserve area (: area and scale distortion grow with the distance from the projection's center point.)} $\square $\\
\end{rem}

Apparently, for our case (and not only), the most interesting property of generalized stereographic projection is its \textit{conformality}. The interested reader can find a proof of this property in the Chapter 3 of the book \cite{GehringMartinPalka} by Frederick W. Gehring, Gaven J. Martin, and Bruce P. Palka. The contents of a significant part of the latter portions of this Chapter are a synthesis of material found in Lars Ahlfors’ Ordway Lectures \cite{Ahlfors} or in Alan Beardon’s book The Geometry of Discrete Groups \cite{Beardon}.An important and noteworthy generalization of all this in the more general framework of a Hilbert space was given in 2018 by Harry Gingold, Yotam Gingold and Salah Hamad in Theorem 1.1 of their paper \cite{GingoldHGingolYdHamad}. However, for the sake of completeness, we will give a sketch of the steps of two alternative approaches to the proof of the conformality of the generalized stereographic projection $\mathcal{P}_{D}: (\mathbb{R}^{D+1}\supset)\hspace{0.1cm} \mathcal{S}^{D} \setminus \lbrace \mathsf{N} \rbrace \longrightarrow \mathbb{R}^{D} (\subset \mathbb{R}^{D+1}) $.

\begin{theorem}
\normalfont{(\textbf{Conformallity of generalized stereographic projection})}. \textit{The generalized stereographic projection is conformal.}
\end{theorem}
\textbf{Proof}. Generalized stereographic projection is conformal because it maps circles on the $D-$sphere to circles and lines on the $D-$plane $\mathbb{R}^{D}$ of $\mathbb{R}^{D+1}$, and it preserves angles locally between intersecting curves on the hypersphere. This angle preservation stems from the property that the mapping is the restriction of a homothety (scaling) and an inversion, which are conformal transformations. More formally, the induced metric on the plane is a scaled version of the pullback metric from the sphere, making the projection conformal.\\

\noindent \underline{First Method}. Sketch of Proof (Intuitive Approach):\textit{ The Geometric Argument }
\begin{itemize}
\item \textit{$1^{st}$ Step: Angle Preservation by Projection.} Consider two intersecting curves on the punctured $D-$sphere $\mathcal{S}^{D} \setminus \lbrace \mathsf{N} \rbrace$ and their corresponding curves after stereographic projection onto the $D-$hyperplane $\mathbb{R}^{D}$ of $\mathbb{R}^{D+1}$. 
\item \textit{$2^{nd}$ Step: Tangent Planes and Radii.} At the intersection point $\Sigma$ on $\mathcal{S}^{D} \setminus \lbrace \mathsf{N} \rbrace$, the tangent plane $\mathcal{T}_{\Sigma} (\mathcal{S}^{D} \setminus \lbrace \mathsf{N} \rbrace)$ and the lines from the projection center $\mathsf{N}$ to $\Sigma$ form certain angles.
\item \textit{$3^{rd}$ Step: Parallel Planes.} When these tangent lines and curves are projected, the resulting tangent lines on  $\Sigma$ and $\mathsf{N}$.
\item \textit{$4^{th}$ Step: Angle Equivalence.} Due to the nature of the projection, the angles between these intersecting curves on $\mathcal{S}^{D} \setminus \lbrace \mathsf{N} \rbrace$ are equal to the angles between their projected counterparts on $\mathbb{R}^{D}$. 
\item \textit{$5^{th}$ Step: Homothety.} The effect of moving the target hyperplane is a homothety (scaling), which is itself a conformal map. 
\end{itemize}

\noindent \underline{Second Method}. Sketch of Proof (Formal Approach):\textit{ Metric-Based Argument }
\begin{itemize}
\item \textit{$1^{st}$ Step: Define the Metrics.} Consider the Euclidean metric $g_{\mathbb{R}^{D}}$ on $\mathbb{R}^{D}$ and the spherical metric $g_{\mathcal{S}^{D}}$ on $\mathcal{S}^{D}$ \footnote{The notation $g_{\mathcal{S}^{D}}$ likely refers to a Riemannian metric on the $D-$sphere $\mathcal{S}^{D}$, which is a manifold of dimension $D$. The standard way to define such a metric, also known as the canonical spherical metric, is to pull back the standard Euclidean metric from the ambient Euclidean space (like $\mathbb{R}^{D+1}$) onto the hypersphere via the inclusion map. This metric gives the hypersphere the same geometry as the unit $D-$sphere in $\mathbb{R}^{D+1}$, and is used to define the spherical measure and the spherical distance. }.
 
\item \textit{$2^{nd}$ Step: Inverse Projection.} The inverse stereographic projection, $\mathcal{P}_{D}^{-1}$, is used to define the metric on the hypersphere.

\item \textit{$3^{rd}$ Step: Pullback Metric.} The conformal nature of the projection can be seen by computing the pullback of the Euclidean metric on $\mathbb{R}^{D}$ to $\mathcal{S}^{D} \setminus \lbrace \mathsf{N} \rbrace$. 

\item \textit{$4^{th}$ Step: Induced Metric.} The inverse map $\mathcal{P}_{D}^{-1}$ allows us to compare the metric on the plane to the metric on the hypersphere. 
 
\item \textit{$5^{th}$ Step: Conformal Factor.} For a point $\Sigma^{'}\in \mathbb{R}^{D}$ and a tangent vector $V\in\mathcal{T}_{\Sigma^{'}} (\mathbb{R}^{D})$, the pullback metric is given by 
\begin{center}
 $(\mathcal{P}_{D}^{-1})_{g_{\mathcal{S}^{D}}}^{*}(V,V)$.
\end{center}

\item \textit{$6^{th}$ Step: Relationship.} It can be shown that this pullback metric on $\mathcal{S}^{D} \setminus \lbrace \mathsf{N} \rbrace$ is related to the Euclidean metric $g_{\mathbb{R}^{D}}$ by a scaling factor, specifically  
\begin{center}
 $g_{\mathcal{S}^{D}}=\frac{4}{(\parallel u \parallel^{2}+1)^{2}}(\mathcal{P}_{D}^{-1})_{g_{\mathcal{S}^{D}}}^{*}$,
\end{center}
where $(u, 0)$ is the point $\Sigma_{0}^{'}$ where the line through $\mathsf{N}$ to the point $\Sigma \in \mathcal{S}^{D} \setminus \lbrace \mathsf{N} \rbrace$ intersects the hyperplane $\mathbb{R}^{D}=\lbrace x_{D+1} = 0\rbrace$ in $\mathbb{R}^{D+1}$. The scaling factor indicates that the map is conformal, preserving angles but not necessarily lengths \footnote{Since, as already noted, the notation $g_{\mathcal{S}^{D}}$ likely refers to a Riemannian metric on the $D-$ sphere $\mathcal{S}^{D}$, let us briefly see how the algorithmic steps of a constructive definition of the Riemannian metric $g_{S^{D}}$ are structured.
\begin{enumerate}
\item \textit{Define the hypersphere}: The $D-$sphere $\mathcal{S}^{D}$ is typically defined as the set of points in $\mathbb{R}^{D+1}$ that are at a fixed distance (e.g., $1$) from the origin.
\item \textit{Identify the embedding}: The sphere is viewed as a subset of $\mathbb{R}^{D+1}$ using the canonical inclusion map, denoted as $\iota: \mathcal{S}^{D}  \longrightarrow \mathbb{R}^{D+1}$. 
\item \textit{Pull back the Euclidean metric}: The standard Euclidean metric on $\mathbb{R}^{D+1}$ is 
$g_{\mathbb{R}^{D+1}}= \sum_{\nu=1}^{D+1} dx_{\nu}\otimes dx_{\nu}$.
\item \textit{Apply the pullback}: The spherical metric $g_{\mathcal{S}^{D}}$ is then obtained by pulling back the Euclidean metric from $\mathbb{R}^{D+1}$ to $\mathcal{S}^{D}$ via the inclusion map $\iota: \mathcal{S}^{D}  \longrightarrow \mathbb{R}^{D+1}$ as follows. For any two tangent vectors $v$ and $w$ at a point on $\mathcal{S}^{D}$, the dot product is calculated as 
$g_{\mathcal{S}^{D}} (v,w) = g_{\mathbb{R}^{D+1}} (\iota(v),\iota(w))$.
Since the inclusion map $\iota$ is an embedding, its differential $d_{\iota}=\iota^{*}$ is the identity on the tangent spaces ($D_{\iota}$ acts as the identity on tangent vectors, making the pullback the restriction of the Euclidean metric to $\mathcal{S}^{D}$), so $g_{\mathcal{S}^{D}} (v,w)=g_{\mathbb{R}^{D+1}} (v,w)$:
\begin{center}
$g_(\mathcal{S}^{D})(v,w)= g_{\mathbb{R}^{D+1}} (d_{\iota} (v),d_{\iota} (w))=g_{\mathbb{R}^{D+1}} (\iota^{*}(v),\iota^{*}(w))=g_{\mathbb{R}^{D+1}}(v,w)$.
\end{center}
\end{enumerate}
The key properties of the Riemannian metric $g_{\mathcal{S}^{D}}$ are the following.
\begin{itemize}
\item \textit{Induced spherical measure}: The surface measure $\sigma^{D}$ on $\mathcal{S}^{D}$ which represents the natural surface area measure on the $D-$sphere, is induced by this Riemannian metric.
\item \textit{Constant sectional curvature}: The canonical spherical metric is a metric which gives $\mathcal{S}^{D}$ a constant sectional curvature. 
\item \textit{Spherical Coordinates}: In suitable coordinates (like spherical coordinates, or via stereographic projection), the induced metric can be expressed in a concrete form, such as a conformally Euclidean metric in stereographic coordinates. 
\item \textit{This metric is a canonical metric}: This metric is the natural, or canonical, metric on the hypersphere. 
\item \textit{This metric defines spherical distance}: This metric defines the shortest path between two points on the hypersphere, which is the angle between the vectors pointing from the origin to those points.
\item \textit{This metric is homogeneous}: The spherical metric is a homogeneous metric, meaning the hypersphere looks the same at every point.
\end{itemize}
}. $\square $\\ 
\end{itemize}

\section{Conformal dimensionality reduction}

Let
\begin{center}
$\mathbb{X}_{D+1}=\lbrace \textbf{x}^{(\nu )}=\left( x_{1}^{(\nu )}, x_{2}^{(\nu )},...,x_{D}^{(\nu )},x_{D+1}^{(\nu )} \right) \in \mathbb{R}^{D+1}, \nu = 1, 2,..., n \rbrace \left(\subset \mathbb{R}^{D+1} \right)$ 
\end{center}
be any given dataset of $n$ points in $\mathbb{R}^{D+1}$. Based on what was mentioned in the previous Section, we will construct a conformal homeomorphism method for dimensionality reduction.

\begin{theorem}
Let $(2\leq) d < D$. There is a polynomial time dimensionality reduction algorithm transforming any dataset $\mathbb{X}_{D+1}$   into a new isopleth dataset 
\begin{center}
$\mathbb{Y}_{d+1}=\lbrace \textbf{y}^{(\nu )}=\left( y_{1}^{(\nu )}, y_{2}^{(\nu )},...,y_{d}^{(\nu )},y_{d+1}^{(\nu )} \right) \in \mathbb{R}^{d+1}, \nu = 1, 2,..., n \rbrace \left(\subset \mathbb{R}^{d+1} \right)$ 
\end{center} 
in such a way that
\begin{itemize}
\item all the angles of the new shape that will result from the $n$ points of the set $\mathbb{Y}_{d+1}$ within the Euclidean space $\mathbb{R}^{d+1}$ of the smaller dimension $d+1$ keep unchanged all the angles of the original shape that was formed by the $n$ points of the set $\mathbb{X}_{D+1}$ within the Euclidean space $\mathbb{R}^{D+1}$ of the larger dimension $D+1$ and 
\item completely ignoring whether or not there is a minimum number of variables required to represent the data set $\mathbb{X}_{D+1}$ (and thus completely ignoring whether or not there is an intrinsic dimension for the data set $\mathbb{X}_{D+1}$).
\end{itemize}
\end{theorem}
\textbf{Proof}. The construction of the algorithm is done according to the following three Steps.\\

\noindent   \underline{$1^{st}$ Step}. For $\nu=1,2,…,n$ do
\begin{itemize}
\item Put
\begin{center}
$\textbf{x}^{(D+1,\nu)}:=
\left( 
x_{1}^{(D+1,\nu)},...,x_{D}^{(D+1,\nu)},x_{D+1}^{(D+1,\nu)}
\right)
\equiv
\left( 
x_{1}^{(\nu)}, ...,x_{D}^{(\nu)},x_{D+1}^{(\nu)}
\right) \in \mathbb{X}_{D+1} $.
\end{center}
\item Normalize $\textbf{x}^{(D+1,\nu)}$ to obtain $\chi^{(D+1,\nu)}\in \mathcal{S}^{D} \setminus \lbrace \mathsf{N} \rbrace$: 
\begin{center}
$\chi^{(D+1,\nu)}:=\frac{\textbf{x}^{(D+1,\nu)}}{\parallel \textbf{x}^{(D+1,\nu)} \parallel}=
\biggl(
\frac{x_{1}^{(D+1,\nu)}}{\underbrace{\parallel \textbf{x}^{(D+1,\nu)} \parallel}_{:=\chi_{1}^{(D+1,\nu)}}},...,
\frac{x_{D}^{(D+1,\nu)}}{\underbrace{\parallel \textbf{x}^{(D+1,\nu)} \parallel}_{:=\chi_{D}^{(D+1,\nu)}}},
\frac{x_{D+1}^{(D+1,\nu)}}{\underbrace{\parallel \textbf{x}^{(D+1,\nu)} \parallel}_{:=\chi_{D+1}^{(D+1,\nu)}}}
\biggr) \in \mathcal{S}^{D} \setminus \lbrace \mathsf{N} \rbrace
$,
\end{center}
where we have used the notation 
\begin{center}
$\parallel \textbf{x}^{(D+1,\nu)} \parallel:=\sqrt{\sum_{j=1}^{D+1}\mid x_{j}^{(D+1,\nu)} \mid^{2}}$.
\end{center}
\item Map $\chi^{(D+1,\nu)}= \left(  \chi_{1}^{(D+1,\nu)},...,\chi_{D}^{(D+1,\nu)}, \chi_{D+1}^{(D+1,\nu)} \right) $ to 
\begin{center}
$\mathcal{P}_{D}\left(\chi^{(D+1,\nu)}\right) 
:=
\biggl(
\frac{\chi_{1}^{(D+1,\nu)}}{\underbrace{1-\chi_{D+1}^{(D+1,\nu)}}_{:=x_{1}^{(D,\nu)}}},
\frac{\chi_{2}^{(D+1,\nu)}}{\underbrace{1-\chi_{D+1}^{(D+1,\nu)}}_{:=x_{2}^{(D,\nu)}}},...,
\frac{\chi_{D}^{(D+1,\nu)}}{\underbrace{1-\chi_{D+1}^{(D+1,\nu)}}_{:=x_{D}^{(D,\nu)}}}
\biggr)$.
\end{center}
\item Put $\textbf{x}^{(D,\nu)}:= \left(x_{1}^{(D,\nu)},x_{2}^{(D,\nu)},...,x_{D}^{(D,\nu)} \right) \in \mathbb{R}^{D}$.
\end{itemize}

\newpage

\noindent   \underline{$2^{nd}$ Step}. For $\ell=D,D-1,…,d+2$ do\\
\indent \hspace{0.6cm} For $\nu=1,2,…,n$ do\\
\begin{itemize}
\item Normalize $\textbf{x}^{(\ell,\nu)}$ to obtain $\chi^{(\ell,\nu)}\in \mathcal{S}^{\ell-1} \setminus \lbrace \mathsf{N} \rbrace$: 
\begin{center}
$\chi^{(\ell,\nu)}:=\frac{\textbf{x}^{(\ell,\nu)}}{\parallel \textbf{x}^{(\ell,\nu)} \parallel}=
\biggl(
\frac{x_{1}^{(\ell,\nu)}}{\underbrace{\parallel \textbf{x}^{(\ell,\nu)} \parallel}_{:=\chi_{1}^{(\ell,\nu)}}},...,
\frac{x_{\ell-1}^{(\ell,\nu)}}{\underbrace{\parallel \textbf{x}^{(\ell,\nu)} \parallel}_{:=\chi_{\ell-1}^{(\ell,\nu)}}},
\frac{x_{\ell}^{(\ell,\nu)}}{\underbrace{\parallel \textbf{x}^{(\ell,\nu)} \parallel}_{:=\chi_{\ell}^{(\ell,\nu)}}}
\biggr) \in \mathcal{S}^{\ell-1} \setminus \lbrace \mathsf{N} \rbrace
$,
\end{center}
where we have used the notation 
\begin{center}
$\parallel \textbf{x}^{(\ell,\nu)} \parallel:=\sqrt{\sum_{j=1}^{\ell}\mid x_{j}^{(\ell,\nu)} \mid^{2}}$.
\end{center}
\item Map $\chi^{(\ell,\nu)}= \left(  \chi_{1}^{(\ell,\nu)},...,\chi_{\ell-1}^{(\ell,\nu)}, \chi_{\ell}^{(\ell,\nu)} \right) $ to 
\begin{center}
$\mathcal{P}_{\ell-1}\left(\chi^{(\ell,\nu)}\right) 
:=
\biggl(
\frac{\chi_{1}^{(\ell,\nu)}}{\underbrace{1-\chi_{\ell}^{(\ell,\nu)}}_{:=x_{1}^{(\ell-1,\nu)}}},
\frac{\chi_{2}^{(\ell,\nu)}}{\underbrace{1-\chi_{\ell}^{(\ell,\nu)}}_{:=x_{2}^{(\ell-1,\nu)}}},...,
\frac{\chi_{\ell-1}^{(\ell,\nu)}}{\underbrace{1-\chi_{\ell}^{(\ell,\nu)}}_{:=x_{\ell-1}^{(\ell-1,\nu)}}}
\biggr) \in \mathbb{R}^{\ell-1}$.
\end{center}
\item Put $\textbf{x}^{(\ell-1,\nu)}:= \left(x_{1}^{(\ell-1,\nu)},x_{2}^{(\ell-1,\nu)},...,x_{\ell-1}^{(\ell-1,\nu)} \right) \in \mathbb{R}^{\ell-1}$. 
\end{itemize}

\noindent   \underline{$3^{rd}$ Step}. For $\nu=1,2,…,n$ put
\begin{center}
$\textbf{y}^{(\nu )}:= \left( y_{1}^{(\nu )},...,y_{d}^{(\nu )}, y_{d+1}^{(\nu )} \right)\equiv 
\left(x_{1}^{(d+1,\nu)},x_{2}^{(d+1,\nu)},...,x_{d+1}^{(\ell-1,\nu)} \right) \in \mathbb{R}^{d+1}$.
\end{center}

The new dataset
\begin{center}
$\mathbb{Y}_{d+1}=\lbrace \textbf{y}^{(\nu )}=\left( y_{1}^{(\nu )}, y_{2}^{(\nu )},...,y_{d}^{(\nu )},y_{d+1}^{(\nu )} \right) \in \mathbb{R}^{d+1}, \nu = 1, 2,..., n \rbrace \left(\subset \mathbb{R}^{d+1} \right)$ 
\end{center}
is the final and desired set, since all the angles of the new shape that will result from the $n$ points of the set $\mathbb{Y}_{d+1}$ within the Euclidean space $\mathbb{R}^{d+1}$ of the smaller dimension $d+1$ keep unchanged all the angles of the original shape that was formed by the $n$ points of the set $\mathbb{X}_{D+1}$ within the Euclidean space $\mathbb{R}^{D+1}$ of the larger dimension $D+1$. Notice that this method completely ignores whether or not there is a minimum number of variables required to represent the dataset $\mathbb{X}_{D+1}$ (and thus completely ignores whether or not there is an intrinsic dimension for the dataset $\mathbb{X}_{D+1}$.\\

Regarding the computational complexity, it is clear that, in the first Step, $n(D+1)$ multiplications, $n(D+1)$ additions, $n$ subtractions, $n$ square root extractions and finally $n(D+1)$ divisions are required, accompanied by $nD$ divisions. That is, in the first Step a total of  $4nD+5n$ operations are required. Similarly, in each iteration $\ell=D, D-1, ..., d+2$ of the second Step,  $4n\ell+5n$ operations are required. In total, in the first Step and in all iterations of the second Step,\\

\noindent \hspace{0.6cm} $\left( 4D+5+4\left( D-1\right) +5+4\left( D-2\right) +5+...+4\left( D-\left( d+1\right) \right) +5\right) n$\\

\noindent \hspace{7.6cm} $=\left( \sum_{\ell=0}^{D-(d+1)}\left(4\left( D-\ell \right) +5 \right) \right)n= \textbf{O}\left(n D^{2} \right)  $\\

\noindent operations are required. This shows the polynomial complexity of the method. $\square $\\

\section{Conformal dimensionality increase}

Dimensionality increase refers to a dataset having more features or dimensions, which is often associated with the ``\textit{curse of dimensionality}" (\cite{Bellman}) - a phenomenon where exponential growth in the data's volume leads to increased data sparsity, similar distances between points, higher computational costs, and a greater risk of model overfitting. While higher dimensionality can capture more complex data structures, it poses significant challenges for pattern recognition and analysis, requiring specialized techniques or vast amounts of data to maintain performance.\\
 
Even though the high dimensionality of input is usually supposed to make learning of RL agents more difficult, in \cite{OtaOikiJhaMariyamaNikovski} it is showed that the reinforcement learning (RL) agents in fact learn more efficiently with the high-dimensional representation than with the lower-dimensional state observations. The authors Kei Ota, Tomoaki Oiki, Devesh K. Jha, Toshisada Mariyama and Daniel Nikovski believe that \textit{stronger feature propagation together with larger networks (and thus larger search space) allows RL agents to learn more complex functions of states and thus improves the sample efficiency}. Through numerical experiments, they showed that the proposed method outperforms several other state-of-the-art algorithms in terms of both sample efficiency and performance. Next, in 2024, Kei Ota, Devesh K. Jha and Asako Kanezaki used a three-fold technique, to show that one can train very large networks that result in significant performance gains (\cite{OtaJhaKanezaki}). Regardless of all this, Tatyana Barron argued that, formally, a finite-dimensional mathematical model for finding solutions to applied problems in data spaces can be embedded into a higher-dimensional model inside of which a desired solution will exist (\cite{OtaOikiJhaMariyamaNikovski}).\\

Let
\begin{center}
$\mathbb{X}_{D+1}=\lbrace \textbf{x}^{(\nu )}=\left( x_{1}^{(\nu )}, x_{2}^{(\nu )},...,x_{D}^{(\nu )},x_{D+1}^{(\nu )} \right) \in \mathbb{R}^{D+1}, \nu = 1, 2,..., n \rbrace \left(\subset \mathbb{R}^{D+1} \right)$ 
\end{center}
be any given dataset of $n$ points in $\mathbb{R}^{D+1}$. Based on what was mentioned in Section 1, we will construct a conformal homeomorphism method for dimensionality increase.

\begin{theorem}
Let $D < D^{'}$. There is a polynomial time dimensionality increase algorithm transforming the dataset $\mathbb{X}_{D+1}$   into a new isopleth dataset 
\begin{center}
$\mathbb{Y}_{D^{'}+1}=\lbrace \textbf{y}^{(\nu )}=\left( y_{1}^{(\nu )}, y_{2}^{(\nu )},...,y_{D^{'}}^{(\nu )},y_{D^{'}+1}^{(\nu )} \right) \in \mathbb{R}^{D^{'}+1}, \nu = 1, 2,..., n \rbrace \left(\subset \mathbb{R}^{D^{'}+1} \right)$ 
\end{center} 
in such a way that any angle of the new shape resulting from $n$ points of the set $\mathbb{Y}_{D^{'}+1}$ within the Euclidean space $\mathbb{R}^{D^{'}+1}$ of larger dimension $D^{'}+1$ is equal to the angle of the shape that was formed by the $n$ corresponding archetypal points of the original dataset $\mathbb{X}_{D+1}$ within the Euclidean space $\mathbb{R}^{D+1}$ of smaller dimension $D+1$.
\end{theorem}
\textbf{Proof}. The construction of the algorithm is done according to the following two Steps.\\

\noindent   \underline{$1^{st}$ Step}. For $\ell=D,D+1,…,D^{'}$ do\\
\indent \hspace{0.6cm} For  $\nu=1,2,…,n$ do

\begin{itemize}
\item Put
\begin{center}
$\textbf{x}^{(\ell+1,\nu)}:=
\left( 
x_{1}^{(\ell+1,\nu)},...,x_{\ell}^{(\ell+1,\nu)},x_{\ell+1}^{(\ell+1,\nu)}
\right)
\equiv
\left( 
x_{1}^{(\nu)}, ...,x_{\ell}^{(\nu)},x_{\ell+1}^{(\nu)}
\right) \in \mathbb{X}_{\ell+1} $.
\end{center}

\newpage

\item Map $\textbf{x}^{(\ell+1,\nu)}\in \mathbb{R}^{\ell+1}$ to $\mathcal{P}_{\ell}^{-1}\left( \textbf{x}^{(\ell+1,\nu)} \right)   $ to take the new point:
\begin{center}
$\textbf{x}^{(\ell+2,\nu)}
:=
\biggl(
\frac{2 x_{1}^{(\ell+1,\nu)}}{\underbrace{\parallel \textbf{x}^{(\ell+1,\nu)} \parallel^{2}+1 }_{:=x_{1}^{\ell+2,\nu}}},...
\frac{2 x_{\ell}^{(\ell+1,\nu)}}{\underbrace{\parallel \textbf{x}^{(\ell+1,\nu)} \parallel^{2}+1 }_{:=x_{\ell+1}^{\ell+2,\nu}}},
\frac{\parallel \textbf{x}^{(\ell+1,\nu)} \parallel^{2}-1 }{\underbrace{\parallel \textbf{x}^{(\ell+1,\nu)} \parallel^{2}+1 }_{:=x_{\ell+2}^{\ell+2,\nu}}},
\biggr)
\in \mathcal{S}^{\ell+1} \setminus \lbrace \mathsf{N} \rbrace \left( \subset \mathbb{R}^{\ell+2}\right) $.
\end{center}
where we have used the notation 
\begin{center}
$\parallel \textbf{x}^{(\ell+1,\nu)} \parallel^{2}:=\sum_{j=1}^{\ell+1}\mid x_{j}^{(\ell+1,\nu)} \mid^{2}$.
\end{center}
\item Put $\textbf{x}^{(D^{'}+1,\nu)}:= \left(x_{1}^{(D^{'}+1,\nu)},...,x_{D^{'}}^{(D^{'}+1,\nu)},x_{D^{'}+1}^{(D^{'}+1,\nu)} \right) \in \mathbb{R}^{D^{'}+1}$.
\end{itemize}
\noindent   \underline{$2^{nd}$ Step}. For $\nu=1,2,…,n$ put
\begin{center}
$\textbf{y}^{(\nu )}:= \left( y_{1}^{(\nu )},...,y_{D^{'}}^{(\nu )}, y_{D^{'}+1}^{(\nu )} \right)\equiv 
\textbf{x}^{(D^{'}+1,\nu)} \in \mathbb{R}^{D^{'}+1}$.
\end{center}
The new dataset
\begin{center}
$\mathbb{Y}_{D^{'}+1}=\lbrace \textbf{y}^{(\nu )}=\left( y_{1}^{(\nu )}, y_{2}^{(\nu )},...,y_{D^{'}}^{(\nu )},y_{D^{'}+1}^{(\nu )} \right) \in \mathbb{R}^{D^{'}+1}, \nu = 1, 2,..., n \rbrace \left(\subset \mathbb{R}^{D^{'}+1} \right)$ 
\end{center}
is the desired set, since all the angles of the new shape that will result from the $n$ points of the set $\mathbb{Y}_{D^{'}+1}$ within the Euclidean space $\mathbb{R}^{D^{'}+1}$ of larger dimension $D^{'}+1$ keep unchanged all the angles of the original shape that was formed by the $n$ points of the set $\mathbb{X}_{D+1}$ within the Euclidean space $\mathbb{R}^{D+1}$ of the smaller dimension $D+1$ \footnote{As in Theorem 4, notice that this method completely ignores whether or not there is a minimum number of variables required to represent the dataset $\mathbb{X}_{D+1}$ (and thus completely ignores whether or not there is an intrinsic dimension for the dataset $\mathbb{X}_{D+1}$.}.\\

Regarding the computational complexity, it is clear that in each iteration $2n(\ell+1)$ multiplications, $n(\ell+2)$ additions, $n$ subtractions and finally $n(\ell+2)$ divisions are required; in total $n(4\ell +8)$ operations. It follows that in all iterations the number of required numerical operations equals \\

\noindent \hspace{0.6cm} $4n \sum_{\ell=D}^{D^{'}}\ell+8n(D^{'}-D)=4n(\mathfrak{d}+1)D+4n\sum_{\kappa=0}^{\mathfrak{d}}\kappa+8n\mathfrak{d} $\\

\noindent \hspace{7cm} $=4n(\mathfrak{d}+1)D+4n\frac{\mathfrak{d}(\mathfrak{d}+1)}{2}+8n\mathfrak{d} $\\

\noindent \hspace{7cm} $=4n(\mathfrak{d}+1)D+2n\mathfrak{d}(\mathfrak{d}+1)+8n\mathfrak{d} = \textbf{O}\left(n (D^{'}-D)^{2} D \right)$,\\

\noindent where we have used the notation $\mathfrak{d}:=D^{'}-D$. This shows that the complexity of the method is polynomial and the proof is complete. $\square $\\

\noindent
\vspace{-\baselineskip}
\end{document}